\documentclass [12pt]  {article} 
\usepackage{amssymb}
\def \ZZ{{\mathbb{Z}}}
\def \QQ{{\mathbb{Q}}}

\def \FF{{\mathbb{F}}}
\def \PP{{\mathbb{P}}}
\def \TTT{{\cal T}}
\def \tn{{\hbox{$\not|\,$}}}
\def \nn{{\mathfrak n}}
\def \pp{{\mathfrak p}}

\begin{document}

\begin{center}
{\Large {\bf Strong Weil curves over $\FF_q(T)$ with
small conductor}}\\
\bigskip
{\sc Andreas Schweizer}\\
\bigskip
{\small {\rm Institute of Mathematics, Academia Sinica\\
6F, Astronomy-Mathematics Building\\
No. 1, Sec. 4, Roosevelt Road\\
Taipei 10617, Taiwan\\
e-mail: schweizer@math.sinica.edu.tw}}
\end{center}
\begin{abstract}
We continue work of Gekeler and others on elliptic curves
over $\FF_q(T)$ with conductor $\infty\cdot\nn$ where
$\nn\in\FF_q[T]$ has degree $3$. Because of the Frobenius
isogeny there are infinitely many curves in each isogeny 
class, and we discuss in particular which of these curves
is the strong Weil curve with respect to the uniformization
by the Drinfeld modular curve $X_0(\nn)$. As a corollary we 
obtain that the strong Weil curve $E/\FF_q(T)$ always gives 
a rational elliptic surface over $\overline{\FF_q}$.
\end{abstract}

\subsection*{0. Introduction}

The first paper that systematically used Drinfeld modular curves 
and the Bruhat-Tits tree in order to classify elliptic curves 
over a function field was [Ge1]. 
\par
With an eye towards feasibility of practical calculations, what 
we are talking about here are elliptic curves over a rational 
function field $\FF_q(T)$ and Drinfeld modular curves $X_0(\nn)$
for $\nn\in\FF_q[T]$. Every elliptic curve over $\FF_q(T)$ that
is modular, that is, covered by some $X_0(\nn)$, must have split 
multiplicative reduction at the place $\infty$ ($=$ pole divisor
of $T$). Conversely, every elliptic curve over $\FF_q(T)$ with
conductor $\infty\cdot\nn$ and split multiplicative reduction
at $\infty$ is an isogeny factor of the Jacobian of $X_0(\nn)$.
\par
It is known that there are no such elliptic curves if $deg(\nn)\le 2$. 
So elliptic curves with conductor $\infty\cdot\nn$ where $deg(\nn)=3$ 
represent the case with the smallest possible conductor. When considered 
as elliptic surfaces over the algebraic closure of $\FF_q$ they give 
so-called extremal elliptic surfaces, which means that they are 
supersingular (in the sense of surface theory) and their Mordell-Weil 
groups are finite even over $\overline{\FF_q}(T)$ (compare [Shi]).
\par
Now let $\TTT$ be the Bruhat-Tits tree of $GL_2(\FF_q((\frac{1}{T})))$.
Then the homology of the quotient graph $\Gamma_0(\nn)\!\setminus\!\TTT$ 
encodes the splitting of the Jacobian of $X_0(\nn)$. So it plays the same 
role as the space of cusp forms of weight $2$ for $\Gamma_0(N)$ in the 
theory of elliptic curves of conductor $N$ over $\QQ$.
\par
In [Ge1] the structure of the quotient graph 
$\Gamma_0(\nn)\!\setminus\!\TTT$ is determined for $\deg(\nn)=3$.
Then the splitting of the Jacobian is calculated for $q\le 16$.
For elliptic curves with certain conductors explicit equations are 
given uniformly in $q$. In [Lei] these calculations were extended 
and further uniform explicit equations of elliptic curves were given.
\par
Other papers have used Drinfeld modular curves to explicitly determine 
all elliptic curves over $\FF_q(T)$ with conductors $\infty\cdot\nn$ 
of high degree but only few places of bad reduction: [Ge4] for $q$ 
a power of $2$ and $\nn=T^n$; [Ge5] ditto for $q$ a power of $3$; 
and [Sch5] for $q$ even and $\nn=T^n(T-1)$. 
\par 
However, although one of the subjects of [Ge2] is the characterization 
of the strong Weil curve, up to now the only explicit calculations 
of strong Weil curves seem to be three examples in [Ge2], the case 
$q=2$, $deg(\nn)=4$ treated in Chapter 4 of [Sch1] (the results are 
also listed in [Sch2]), and unpublished computer calculations by Udo 
Nonnengardt from around 1995 for $q=2$, $\nn=T^n$ with $n\le 10$.
\par
In contrast to the classical situation of modular elliptic curves 
over $\QQ$, each $\FF_q(T)$-isogeny class contains infinitely many 
non-isomorphic curves. More precisely, they are obtained from finitely 
many curves by repeated application of the Frobenius isogeny.
Papikian [Pa1], [Pa2] has shown that in certain situations the strong 
Weil curve is not the Frobenius of another curve over $\FF_q(T)$, but 
from examples it is known that this is not a general phenomenon.
\par
The Frobenius isogeny is more problematic than the other ones,
for example in the following context. By the $2$-dimensional 
L\"uroth theorem, the image of a rational elliptic surface under
a separable isogeny is again a rational elliptic surface. But the
Frobenius of a rational elliptic surface need not be a rational 
surface. For example, the Frobenius of a semistable rational 
elliptic surface is never a rational surface (although it is of 
course a unirational surface).
\par
Now Frobenius minimal extremal elliptic surfaces were proved to 
be rational surfaces in [Ito] for characteristic $p\ge 5$ and in 
[Sch6] for characteristics $2$ and $3$. Moreover, extremal rational 
elliptic surfaces have been explicitly classified in [La1] and [La2]. 
This invites of course the question where the strong Weil curves 
of conductor $\infty\cdot\nn$ with $deg(\nn)=3$ are standing in 
this respect.
\par
This question will be answered quite explicitly in this paper.
The main reason why this can be done for all $q$, not just for 
specific instances, is that for $deg(\nn)=3$ the corresponding
quotient graph $\Gamma_0(\nn)\!\setminus\!\TTT$ can be described
uniformly in $q$ and the necessary calculations and arguments 
can be carried out more or less generically.
\\

\subsection*{1. Basic facts}

For an elliptic curve $E$ over any field $K$ of characteristic 
$p$ we write $E^{(p)}$ for the image of $E$ under the Frobenius 
isogeny. Since $E^{(p)}$ is obtained by raising the coefficients 
in a Weierstra\ss\ equation of $E$ to the $p$-th power, we have 
$j(E^{(p)})=(j(E))^p$. Conversely, if $j(E)\in(K^*)^p$, then 
$E\cong\widetilde{E}^{(p)}$ for some $\widetilde{E}$ over $K$.
\par
If $p\ge 5$, then $Y^2=X^3+a_4X+a_6$ is isomorphic
over $K$ to 
$$Y^2=X^3+a_4^{p^2}X+a_4^{\frac{3(p^2-1)}{2}}a_6$$
and from the formula for the $j$-invariant we see that 
if $j(E)$ is a $p$-th power, then 
$a_4^{\frac{3(p^2-1)}{2}}a_6$ also is.
\par
Analogous arguments work in characteristic $3$, where $E$ has 
a normal form $Y^2=X^3+a_2X^2-\frac{a_2^3}{j(E)}$ with $a_2\neq 0$ 
(see [Si, Appendix A]) and in characteristic $2$, where
$Y^2+XY=X^3+a_2X^2+\frac{1}{j(E)}$ is isomorphic over
$K$ to $Y^2+XY=X^3+a_2^2X^2+\frac{1}{j(E)}$ 
(also [Si, Appendix A]).
\par
Somewhat similar arguments show that an elliptic curve with
$j$-invariant $0$ always is the Frobenius of some elliptic
curve over the same field.
\par
We are in particular interested in elliptic curves over the 
rational function field $\FF_q(T)$ where $\FF_q$ is the finite
field with $q$ elements and characteristic $p$. By the above,
such a curve $E$ over $\FF_q(T)$ with $j(E)\not\in\FF_q$ is
{\bf Frobenius minimal}, i.e. not the Frobenius of another 
elliptic curve over $\FF_q(T)$, if and only if $j(E)$ is not 
a $p$-th power in $\FF_q(T)$.
\\ \\
Now let us fix the polynomial ring $\FF_q[T]$ and hence the place
$\infty$ (the pole divisor of $T$). Then it is known that the 
elliptic curves over $\FF_q(T)$ that are modular, i.e. that are
images of a Drinfeld modular curve $X_0(\nn)$ under a nonconstant 
$\FF_q(T)$-rational morphism, are exactly the ones that have split 
multiplicative reduction at $\infty$.
See [GeRe] for a thorough treatment (over any congruence function 
field), and [Ge2] or Ge3] for less technical accounts over $\FF_q(T)$.
\par
As in the classical situation of elliptic curves over $\QQ$,
in every $\FF_q(T)$-isogeny class of modular elliptic curves of 
conductor $\infty\cdot\nn$ there exists a unique curve $E$ that 
is ``closest to $X_0(\nn)$'' in the sense that every morphism
$X_0(\nn)\to E'$ where $E'$ is in the given isogeny class factors
over $E$ ([GeRe, section 8.4]). Equivalently, $E$ is a subvariety 
of the Jacobian of $X_0(\nn)$, not just an isogeny factor. To see
this, note that since the map $X_0(\nn)\to E$ does not factor over
an elliptic curve, Picard functoriality induces an embedding of 
$E\cong Jac(E)$ into $Jac(X_0(\nn))$. For the converse we are using 
that $Jac(X_0(\nn))$ contains only one abelian subvariety isogenous 
to $E$. This elliptic curve $E$ is called the {\bf strong Weil curve} 
(or optimal elliptic curve in [Pa1]). 
\par
Note however that the notion of strong Weil curve and the degree
of the strong Weil uniformization $\pi: X_0(\nn)\to E$ also depend
on $q$. The base change of the strong Weil curve $E$ over $\FF_q(T)$ 
to $\FF_{q^n}(T)$ is not necessarily the strong Weil curve in its 
$\FF_{q^n}(T)$-isogeny class; and even if it is, the degree of the 
strong Weil uniformization might not be the same. This comes from 
the simple fact that every $q$ has its own Drinfeld modular curves. 
No useful relation is known between the curves $X_0(\nn)$ for the 
same $\nn$ but different $q$ of the same characteristic. 
\par
Naively, one would perhaps guess that the strong Weil curve is always 
Frobenius minimal. This is true for example if $\nn$ is irreducible 
[Pa1, Theorem 1.2] but not in general [GeRe, Examples 9.7.2 and 9.7.3] 
or equivalently [Ge2, Examples 4.2 and 4.3].
\\ \\
In practice, strong Weil curves over $\FF_q(T)$ with conductor 
$\infty\cdot\nn$ can be determined by diagonalizing Hecke operators 
on the homology $H_1(\Gamma_0(\nn)\!\setminus\!\TTT,\ZZ)$ of the
(essentially finite) quotient graph $\Gamma_0(\nn)\!\setminus\!\TTT$.
Here $\TTT$ is the Bruhat-Tits tree of $GL_2(\FF_q((\frac{1}{T})))$ 
and $\Gamma_0(\nn)$ is a Hecke subgroup of $GL_2(\FF_q[T])$.
One obtains a bijection between their $\FF_q(T)$-isogeny classes
(and hence their strong Weil curves) and certain one-dimensional
simultaneous eigenspaces in $H_1(\Gamma_0(\nn)\!\setminus\!\TTT,\ZZ)$.
\par
We suppress discussing the somewhat involved general details [GeRe],
[Ge2], [Ge3] and concentrate on the case of interest to us, that is,
the case $deg(\nn)=3$. 
Up to affine transformation $T\mapsto aT+b$ with $a\in\FF_q^*$, 
$b\in\FF_q$, there are $5$ different cases, namely
\par
$\nn=T^3$ (In this case elliptic curves exist only in characteristics
$2$ and $3$);
\par
$\nn=T^2(T-1)$;
\\
and the three semistable cases
\par
$\nn=T(T-1)(T-c)$, \ $\nn=T\pp_2$, \ $\nn=\pp_3$
\\
where $1\neq c\in\FF_q^*$ and $\pp_i\in\FF_q[T]$ is irreducible 
of degree $i$. 
\par
The corresponding quotient graphs have already been calculated in [Ge1].
For example, the graph $\Gamma_0(T^2(T-1))\!\setminus\!\TTT$ looks 
as follows.
\\
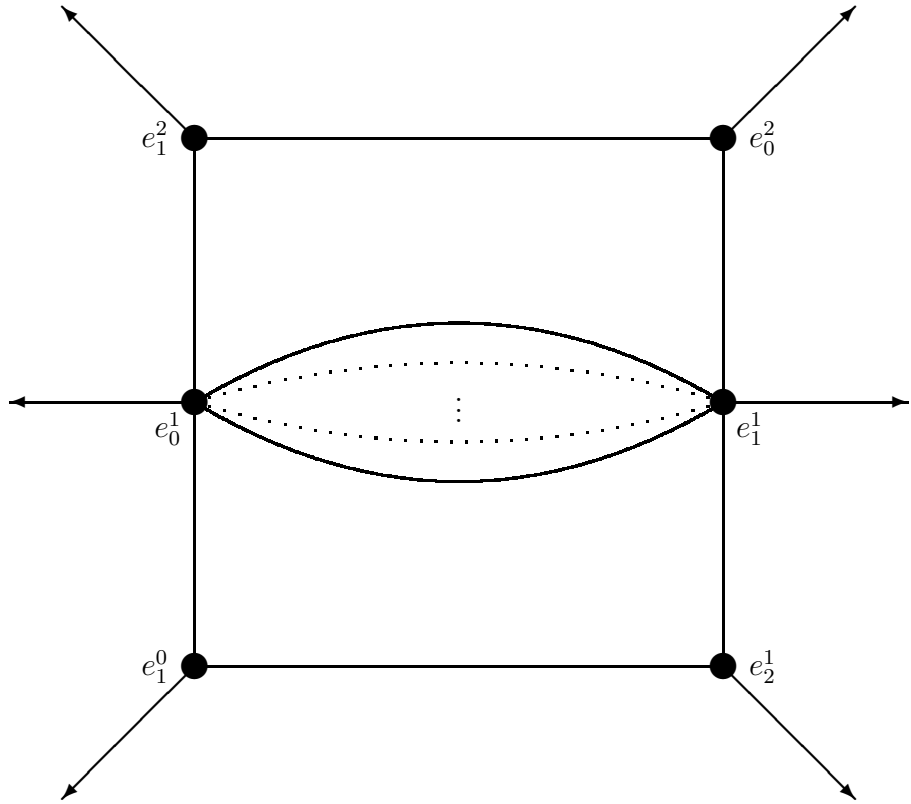
\begin{center}
\setlength{\unitlength}{1pt} \thicklines
\begin{picture}(0,300)
\put(-100,50){\circle*{10}}
\put(-100,150){\circle*{10}}
\put(-100,250){\circle*{10}}
\put(100,50){\circle*{10}}
\put(100,150){\circle*{10}}
\put(100,250){\circle*{10}}
\put(-100,50){\line(1,0){200}} 
\put(-100,250){\line(1,0){200}} 
\put(-100,50){\line(0,1){200}}
\put(100,50){\line(0,1){200}}
\put(-100,250){\vector(-1,1){50}}
\put(-100,50){\vector(-1,-1){50}}
\put(100,250){\vector(1,1){50}}
\put(100,50){\vector(1,-1){50}}
\put(-100,150){\vector(-1,0){70}}
\put(100,150){\vector(1,0){70}}
\qbezier(-100,150)(0,90)(100,150)
\qbezier(-100,150)(0,210)(100,150)
\qbezier[35](-100,150)(0,120)(100,150)
\qbezier[35](-100,150)(0,180)(100,150)
\put(-115,50){\makebox(0,0){$e_1^0$}}
\put(-115,250){\makebox(0,0){$e_1^2$}}
\put(115,50){\makebox(0,0){$e_2^1$}}
\put(115,250){\makebox(0,0){$e_0^2$}}
\put(-110,140){\makebox(0,0){$e_0^1$}}
\put(110,140){\makebox(0,0){$e_1^1$}}
\put(0,150){\makebox(0,0){$\vdots$}}
\end{picture}
\end{center}
Here the dotted lines mean that there are $q-2$ different edges 
between the vertices $e_0^1$ and $e_1^1$. The (in English) somewhat 
unfortunate choice of the letter $e$ for the labeling of the vertices
originates from [Ge1] being written in German.
The $6$ arrows in the picture represent so-called cusps, that is 
half-lines consisting of infinitely many edges. 
\par
For general $\nn\in\FF_q[T]$ of degree $3$ the graph 
$\Gamma_0(\nn)\!\setminus\!\TTT$ looks somewhat similar and can 
be precisely described as follows: There are two vertices $e_0^1$
and $e_1^1$ of valency $q+1$. So the edges attached to $e_0^1$ can
be parametrized by $\PP^1(\FF_q)$. Of these edges the ones that
directly connect $e_0^1$ and $e_1^1$ are those corresponding to the
$a\in\FF_q$ with $(T-a)\tn\nn$. If $\nn$ has a multiple zero, there
is a cusp attached to each of $e_0^1$ and $e_1^1$, otherwise not.
Finally, for each simple zero of $\nn$ in $\FF_q$ and for $\infty$
there is a three edge path between $e_0^1$ and $e_1^1$ with two
cusps attached to it (as in the above picture the upper and the 
lower part of the graph, belonging to $0$ respectively $\infty$.)
\par
So for example, if $\nn\in\FF_q[T]$ is irreducible of degree $3$,
the quotient graph $\Gamma_0(\nn)\!\setminus\!\TTT$ is
\\
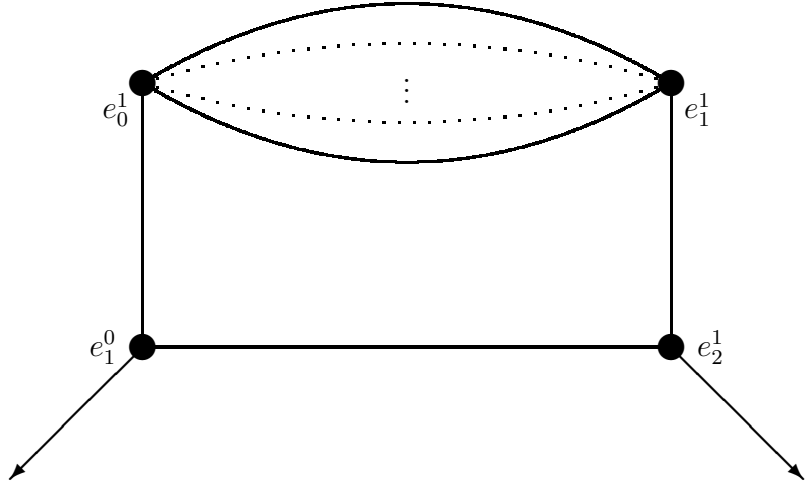
\begin{center}
\setlength{\unitlength}{1pt} \thicklines
\begin{picture}(0,170)
\put(-100,50){\circle*{10}}
\put(-100,150){\circle*{10}}
\put(100,50){\circle*{10}}
\put(100,150){\circle*{10}}
\put(-100,50){\line(1,0){200}} 
\put(-100,50){\line(0,1){100}}
\put(100,50){\line(0,1){100}}
\put(-100,50){\vector(-1,-1){50}}
\put(100,50){\vector(1,-1){50}}
\qbezier(-100,150)(0,90)(100,150)
\qbezier(-100,150)(0,210)(100,150)
\qbezier[35](-100,150)(0,120)(100,150)
\qbezier[35](-100,150)(0,180)(100,150)
\put(-115,50){\makebox(0,0){$e_1^0$}}
\put(115,50){\makebox(0,0){$e_2^1$}}
\put(-110,140){\makebox(0,0){$e_0^1$}}
\put(110,140){\makebox(0,0){$e_1^1$}}
\put(0,150){\makebox(0,0){$\vdots$}}
\end{picture}
\end{center}
where the dotted lines now represent $q$ different edges 
in total between $e_0^1$ and $e_1^1$.
\par
Although we don't need it in the sequel we mention that in our 
description of the graphs the full Atkin-Lehner involution is the 
reflection at the vertical axis whereas in [Ge1] they are printed 
in a way that shows how they project onto the quotient graph 
$GL_2(\FF_q[T])\!\setminus\!\TTT$. Note however that the 
Figures 3 and 5 in [Ge1] have unfortunately been interchanged.
\par
In particular, for $deg(\nn)=3$ the structure of the graph
$\Gamma_0(\nn)\!\setminus\!\TTT$ depends only on the splitting 
type of $\nn$ and can be uniformly described in $q$, whereas
matters are (much) more complicated for $deg(\nn)\ge 4$ [GeNo].
\par
The homology $H_1(\Gamma_0(\nn)\!\setminus\!\TTT,\ZZ)$ of the 
graph $\Gamma_0(\nn)\!\setminus\!\TTT$ is the $\ZZ$-module of 
all $\ZZ$-valued functions $\psi$ on the oriented edges of 
$\Gamma_0(\nn)\!\setminus\!\TTT$ that satisfy the following 
three conditions:
\begin{itemize}
\item[(i)] $\psi(\overline{e})=-\psi(e)$ where $\overline{e}$ 
is the edge $e$ with reversed orientation.
\item[(ii)] (harmonicity) $\sum\limits_{t(e)=v}\psi(e)=0$ for every 
vertex $v$. Here $t(e)$ denotes the terminal vertex of the edge $e$.
\item[(iii)] $\psi$ has finite support. Because of the harmonicity 
this means that $\psi$ vanishes on the cusps.
\end{itemize}
The above description shows that every 
$\psi\in H_1(\Gamma_0(\nn)\!\setminus\!\TTT,\ZZ)$ 
is completely determined by the values it takes on the edges 
terminating in the vertex $e_0^1$. Writing $\psi(a)$ for the 
value of $\psi$ on the edge corresponding to $a\in\PP^1(\FF_q)$ 
establishes a $\ZZ$-module isomorphism between 
$H_1(\Gamma_0(\nn)\!\setminus\!\TTT,\ZZ)$ 
and the $\ZZ$-valued functions $\psi$ on $\PP^1(\FF_q)$ whose
values sum up to $0$ (harmonicity at $e_0^1$), subject to the 
additional condition $\psi(0)=0$ in case $\nn$ has a multiple 
zero, which we place at $0$.
\par
In particular, the dimension of 
$H_1(\Gamma_0(\nn)\!\setminus\!\TTT,\ZZ)$,
which is also the genus of $X_0(\nn)$, 
is $q$ if $\nn$ is square-free, and $q-1$ if not.
\\ \\
{\bf Theorem 1.1.} [Ge1], [Ge2] \it\\
Let $E$ be a strong Weil curve over $\FF_q(T)$ with conductor
$\infty\cdot\nn$ where $deg(\nn)=3$. Let $\varphi$ be the primitive 
cycle in $H_1(\Gamma_0(\nn)\!\setminus\!\TTT,\ZZ)$ belonging to $E$. 
Scaling to $\varphi(\infty)=-1$ we have
$$\varphi(a)=\left\{\begin{array}{ll}
-\lambda_a\ & \hbox{\rm if $E$ has good reduction at $T-a$},\\
-1\  & \hbox{\rm if $E$ has split multiplicative reduction at $T-a$},\\
1\   & \hbox{\rm if $E$ has non-split multiplicative reduction at $T-a$},\\
0\   & \hbox{\rm if $E$ has additive reduction at $T-a$}.
\end{array}\right.$$
(Recall that $\lambda_a =q+1-\#E(\FF_q[T]/(T-a))$ for $(T-a)\tn\nn$.)
Moreover, 
$$-v_{\infty}(j(E))=\min\{\langle\varphi,\psi\rangle>0\ 
:\ \psi\in H_1(\Gamma_0(\nn)\!\setminus\!\TTT,\ZZ)\}$$
where the scalar product is given by
$$\langle\varphi,\psi\rangle=
\sum_{a\in\PP^1(\FF_q)}w_a \varphi(a)\psi(a)$$ 
with $w_a=q+1$ if $E$ has multiplicative reduction at $a\in\PP^1(\FF_q)$, 
and $w_a=1$ otherwise. 
The degree of the strong Weil uniformization $\pi: X_0(\nn)\to E$ is
$$deg(\pi)=\frac{\langle\varphi,\varphi\rangle}{-v_{\infty}(j(E))}.$$
\\ 
{\bf Proof.} \rm \\
The calculation of $\varphi$ is done in [Ge1, Satz 9.1 and p.141].
In particular, $\varphi(a)$ is the negative of the eigenvalue of
$\varphi$ under the Hecke operator $H_{T-a}$. In order to maintain 
compatibility with the examples in [Ge1] and [Lei] we have refrained
from scaling the minus sign away. Also note that the harmonicity of 
$\varphi$ at the vertex $e_0^1$ corresponds to the fact that the 
linear coefficient of the $L$-polynomial of $E$ is $0$. 
\par
For the scalar product see [Ge1, Bemerkung 6.9]. The factor 
$w_a =q+1$ comes from the fact that the middle edge of each 
three edge path between $e_0^1$ and $e_1^1$ carries weight $q-1$. 
Compare [GeRe, 3.2.5] and [GeRe, 4.8]. 
\par
The formulas for $-v_{\infty}(j(E))$ [Ge2, Corollary 3.19] and 
$deg(\pi)$ [Ge2, Corollary 3.20] hold for all $\nn\in\FF_q[T]$,
but in general the scalar product is more complicated and there
is no explicit formula for $\varphi$.
\hfill $\Box$
\\ \\
Theorem 1.1 shows that for $deg(\nn)=3$ the eigenform $\varphi$ 
corresponding to a strong Weil curve $E$, and hence also 
$-v_{\infty}(j(E))$ and the degree of the strong Weil uniformization,
are completely determined by the number of points of the reduction 
at the linear places. We make this explicit.
\\ \\
{\bf Corollary 1.2.} \it
Let $E$ be a strong Weil curve over $\FF_q(T)$ with conductor
$\infty\cdot\nn$ where $deg(\nn)=3$. For every $a\in\FF_q$
such that $(T-a)\tn\nn$ let $\#_{a}$ be the number of 
$\FF_q$-rational points of the reduction modulo $T-a$ of
$E$ (or of any other elliptic curve in the same isogeny class).
Let $N$ be the greatest common divisor of all these $\#_{a}$.
\begin{itemize}
\item[a)] If $E$ has no linear places of non-split multiplicative 
reduction, then the pole order of $j(E)$ at the place $\infty$ is $N$.
\item[b)] If $E$ has at least one linear place of non-split multiplicative 
reduction, then $-v_{\infty}(j(E))=\gcd(N,2q+2)$.
\end{itemize}
\smallskip
\noindent
{\bf Proof. }\rm We write $\delta_a -\delta_{\infty}$ for the cycle 
taking the value $1$ on the edge corresponding to $a$, the value $-1$ 
on the edge corresponding to $\infty$, and the value $0$ on the other 
edges terminating in $e_0^1$.
Obviously
$$\{\delta_a -\delta_{\infty}\ :\ a\in\FF_q\ \hbox{\rm and}\ (T-a)^2\tn\nn\}$$
is a $\ZZ$-basis of $H_1(\Gamma_0(\nn)\!\setminus\!\TTT,\ZZ)$ and thus
$$-v_{\infty}(j(E))=\gcd\{\langle\varphi,\delta_a -\delta_{\infty}\rangle\ 
:\ a\in\FF_q\ \hbox{\rm and}\ (T-a)^2\tn\nn\}.$$
By Theorem 1.1 we have
$$\langle\varphi,\delta_a -\delta_{\infty}\rangle=\left\{\begin{array}{ll}
\#_a \ & \hbox{\rm if $E$ has good reduction at $T-a$,}\\
0    \ & \hbox{\rm if $E$ has split multiplicative reduction at $T-a$,}\\
2q+2 \ & \hbox{\rm if $E$ has non-split multiplicative reduction at $T-a$.}
\end{array}\right.$$
\hfill $\Box$
\\ \\
In particular, if there exists a linear place of good, supersingular
reduction, then the strong Weil curve is Frobenius minimal. The same 
holds if $q$ is odd and $E$ has a linear place of non-split multiplicative
reduction.
\par
For convenient use later on we prove the following, presumably 
well-known facts.
\\ \\
{\bf Lemma 1.3.} \it
\begin{itemize}
\item[a)] An elliptic curve $E$ over a perfect field $k$ of 
characteristic $2$ has a $k$-rational $2$-torsion point if and only 
if it has a model 
$Y^2 +XY=X^3 +a_2 X^2 +a_6$
with $a_6 \neq 0$. Moreover, this model has a $k$-rational $4$-torsion
point if and only if $a_2 =b^2 +b$ for some $b\in k$, that
is, if it can be transformed over $k$ into $Y^2 +XY=X^3 +\frac{1}{j(E)}$.
\item[b)] An elliptic curve $E$ over a perfect field $k$ of 
characteristic $3$ has a $k$-rational $3$-torsion point if and only 
if it has a model 
$Y^2 =X^3 +X^2 +a_6$
with $a_6 \neq 0$.
Furthermore, this model has a $k$-rational $9$-torsion point if and
only if $a_6 =u^3 -u$ with $u\in k\setminus\FF_3$.
\item[c)] An elliptic curve $E$ over a perfect field $k$ of 
characteristic $5$ has a $k$-rational $5$-torsion point if and only 
if it has a model 
$Y^2 =X^3 +3X+a_6$
with $a_6 \neq\pm 1$.
\end{itemize}
\smallskip
\noindent
{\bf Proof. }\rm\\
a) If $E$ has a $2$-torsion point, it is not supersingular. So 
$j(E)\neq 0$, and there exists a change of coordinates over $k$ 
that brings $E$ into the desired form (cf. [Si, Appendix A]). 
With the duplication formula [Si, p.59] one easily checks that 
$(0,\sqrt{a_6})$ is a $2$-torsion point of the above equation,
and that the $4$-torsion points are 
$(\sqrt[4]{a_6}, \sqrt{a_6}+b\sqrt[4]{a_6})$ and
$(\sqrt[4]{a_6}, \sqrt{a_6}+(b+1)\sqrt[4]{a_6})$
with $b^2 +b=a_2$.
\\ \\
b) As before, if $E$ has a $3$-torsion point, it is not supersingular, 
i.e. $j(E)\neq 0$, and by [Si, Appendix A] there exists a model
$Y^2 =X^3 +a_2 X^2 +a_6$. 
With the duplication formula [Si, p.59] one easily verifies that 
the $3$-torsion points of this equation are
$(-\sqrt[3]{a_6}, \pm\sqrt{a_2}\sqrt[3]{a_6})$.
But if $a_2$ is a square in $k$, we can scale $X$ and $Y$ over 
$k$ to get the desired form (with new $a_6$).
\par
After some tedious calculation using the addition formulas 
[Si, p.59] one obtains that in characteristic $3$ the 
triplication formula for the $X$-coordinate of a point 
$P=(x,y)$ on the curve $Y^2 =X^3 +X^2 +a_6$ is
$$x[3]P=\frac{x^9 -a_6 x^3 +a_6^3}{(x^3 +a_6)^2}.$$
So if $\delta$ is the $X$-coordinate of the $3$-torsion point
(and hence $a_6 =-\delta^3$), then $\eta$ is the third power
of the $X$-coordinate of a $9$-torsion point if and only if
$$\eta^3 +\delta^3 \eta-\delta^9=\delta(\eta-\delta^3)^2.$$
Elementary transformations show that this equation is equivalent to
$\delta=\lambda-\lambda^3$ with
$$\lambda=\frac{\delta^2}{\eta+\delta^2 -\delta^3}\ \ \ \hbox{\rm
and }\ \ \eta=\frac{\delta^2}{\lambda}+\delta^3 -\delta^2.$$
This proves that the $X$-coordinate of the $9$-torsion point is
$k$-rational if and only if $a_6 =u^3 -u$ (with $u=\lambda^3$).
Then the $Y$-coordinate is automatically $k$-rational. Otherwise
there would be an element in $Gal(\overline{k}/k)$ that maps the
$9$-torsion point to its inverse. But then it would also map the
($k$-rational!) $3$-torsion point to its inverse, contradiction.
\\ \\
c) The condition $a_6 \neq\pm 1$ is equivalent to $\Delta\neq 0$.
The Hasse invariant of a Weierstra\ss\ equation $Y^2 =X^3 +AX+B$
in characteristic $5$ is $H=2A$. If $H$ is a $4$-th power in $k^*$,
then by [Vo, pp.248/249] the $5$-torsion points are $k$-rational.
Conversely, we have to show that if $(x_0 ,y_0)$ is a $k$-rational
$5$-torsion point, then $\sqrt[4]{H}\in k$; then we can carry out
the desired transformation. 
\par
Now $Q=(x_0 \sqrt{H},y_0\sqrt{H}\sqrt[4]{H})$ 
is a $5$-torsion point of 
$$Y^2 =X^3 +3H^2X+BH\sqrt{H}.$$ 
Its Hasse invariant $H^2$ is a $4$-th power in $k(\sqrt{H})$, so
again by [Vo, pp.248/249] $Q$ is $k(\sqrt{H})$-rational. Since 
$y_0 \neq 0$ this means $\sqrt[4]{H}\in k(\sqrt{H})$. If 
$\sqrt{H}\in k$ we are done. If not, we have 
$\sqrt[4]{H}=u+v\sqrt{H}$ with $u,v\in k$, and after squaring
$0=u^2 +v^2H$, so $\sqrt{H}=\frac{u\sqrt{-1}}{v}\in k$ quite
the same.
\hfill $\Box$
\\

\subsection*{2. The main results}

It is well known that elliptic curves over $\FF_q(T)$ with 
conductor $\infty\cdot T^3$ can only exist in characteristics
$2$ and $3$. If $q$ is a power of $2$ or of $3$, then by 
[Ge5, Corollary 6.4] the Jacobian of the Drinfeld modular 
curve $X_0(T^3)$ is isogenous to a product of $q-1$ elliptic 
curves that are explicitly described in [Ge5]. We prove the 
following refinement.
\\ \\
{\bf Theorem 2.1.} \it
\begin{itemize}
\item[a)] If $q$ is a power of $2$, then the strong Weil
curves over $\FF_q(T)$ with conductor $\infty\cdot T^3$ are
$$Y^2+XY=X^3+\frac{c}{T^4}\ \ \hbox{\it with}\ c\in\FF_q^*.$$
\item[b)] If $q$ is a power of $3$, then the strong Weil 
curves over $\FF_q(T)$ with conductor $\infty\cdot T^3$ are
$$Y^2=X^3+X^2-\frac{c}{T^3}\ \ \hbox{\it with}\ c\in\FF_q^*.$$
\end{itemize}
{\bf Proof. }\rm \\
a) By [Sch7, Theorem 5.1] $Y^2 +XY=X^3 +\frac{c}{T}$ is up to
Frobenius the only curve in its isogeny class. So the pole order
of the $j$-invariant of the strong Weil curve at the place $\infty$ 
must be a power of $2$, and we only have to show that this pole 
order is divisible by $4$ but not by $8$. 
\par   
We use Corollary 1.2. At all places $T-a$ with $a\in\FF_q^*$ the 
reduced curve has an $\FF_q$-rational $4$-torsion point. But it 
is well known (see for example [R\"u, Theorem 1b]) that there 
exists an elliptic curve over $\FF_q$ with non-zero $j$-invariant
that has a $4$-torsion point but no $8$-torsion point over $\FF_q$.
By Lemma 1.3 this curve has an equation $Y^2 +XY=X^3 +a_6$, 
so it occurs as the reduction at one of the places $T-a$.
\\ \\
b) The proof is similar to a). By [Sch6, Proposition 4.3] the 
strong Weil curve must be $Y^2 =X^3+X^2 -\frac{c}{T}$ up to 
Frobenius. So we have to show that the reductions mod $T-a$ with
$a\in\FF_q^*$ all have $3$-torsion points but at least one of
them has no $9$-torsion points over $\FF_q$. This follows again
by combining [R\"u, Theorem 1b] and Lemma 1.3.
\hfill $\Box$
\\ \\
In Chapter 3 of [Lei] four equations of elliptic curves over
$\FF_p(T)$ with conductor $\infty\cdot T^2(T-1)$ are given
uniformly in $p\ge 5$, and by different tricks the curves
are shown to be non-isogenous over $\FF_p(T)$ at least for
the $p$ in certain congruence classes. What seems to have
escaped attention is that by the theory of Tate curves one
can easily see that these $4$ curves are non-isogenous over
any field $\FF_q(T)$ (of characteristic $\ge 5$). 
More precisely:
\\ \\
{\bf Theorem 2.2.} \it
If $char(\FF_q)\ge 5$, then there are $4$ isogeny classes
of elliptic curves over $\FF_q(T)$ with conductor
$\infty\cdot T^2(T-1)$ and split multiplicative reduction
at $\infty$. The Frobenius minimal curves in these classes
are given in the table below, with curves in the same
horizontal box belonging to the same isogeny class. 
The numbers $(mn)$ resp. $(mno^*)$ give the pole orders of 
the $j$-invariant at the places $\infty$ and $T-1$
(and $T$ in the class $E_1$).
\\
\def\abstand{$\vrule width 0pt height 22pt $}
$$ \begin{array}{|r|l|c|} \hline
 no. & equation & j-invariant \\ 
\hline\hline
\abstand E_1:\ \ (222^*) & Y^2=X(X+T)(X+T^2) & 
\frac{2^8(T^2-T+1)^3}{T^2(T-1)^2}  \\
\abstand (114^*) & Y^2=X^3-2T(T-2)X^2+T^4X &
\frac{-2^4(T^2-16T+16)^3}{T^4(T-1)} \\
\abstand (141^*) & Y^2=X^3-2T(T+1)X^2+T^2(T-1)^2X &
\frac{2^4(T^2+14T+1)^3}{T(T-1)^4} \\
\abstand (411^*) & Y^2=X^3+2T(2T-1)X^2+T^2X &
\frac{2^4(16T^2-16T+1)^3}{T(T-1)} \\
 & & \\
\hline
\abstand E_2:\ \ (12) & Y^2=X^3-2T^2X^2+T^3(T-1)X & 
\frac{2^6(T+3)^3}{(T-1)^2} \\
\abstand (21) & Y^2=X^3+4T^2X^2+4T^3X & 
\frac{2^6(4T-3)^3}{T-1} \\
 & & \\
\hline
\abstand E_3:\ \ (13) & 
Y^2=X^3-27T^3(T+8)X+54T^4(T^2-20T-8) & 
\frac{3^3T(T+8)^3}{(T-1)^3} \\
\abstand (31) & Y^2=X^3-3T^3(9T-8)X+2T^4(27T^2-36T+8) &
\frac{3^3T(9T-8)^3}{T-1} \\
 & & \\
\hline
\abstand E_4:\ \ (11) & Y^2=X^3-27T^4X+54T^5(T-2) & 
\frac{2^4\cdot3^3T^2}{T-1} \\
 & & \\
\hline
\end{array} $$ 
\\ \\
{\bf Proof. }\rm
The strategy is to first find all Frobenius minimal elliptic curves 
over $\overline{\FF_p}(T)$ with conductor $\infty\cdot T^2(T-1)$ 
where $\overline{\FF_p}$ is the algebraic closure of $\FF_q$. 
Such a curve is an elliptic surface over $\overline{\FF_p}$; and
since the base curve is a projective line and the conductor has 
degree $4$, it is actually an extremal elliptic surface.
By [Ito, Theorem 3.1] Frobenius minimal extremal elliptic surfaces
are extremal rational elliptic surface, and these have been 
completely classified in [La1] and [La2]. 
Most of the equations in the table we essentially got from [Ito] and 
[Lei]. Note however, that one can apply a M\"obius transformation to 
$\overline{\FF_p}(T)$. This means for example that not all authors 
place the additive fiber at $T=0$, and that the same elliptic surface 
might give rise to several non-isomorphic elliptic curves over 
$\overline{\FF_p}(T)$. For example, there are three pairs of curves 
in the table that are connected by the M\"obius transformation 
$T\mapsto\frac{T}{T-1}$ that fixes $0$ and interchanges $1$ and 
$\infty$. Possibly one has to take an unramified quadratic twist to 
make the multiplicative reduction at $\infty$ split. It turns out 
that for each curve one can find an equation already over $\FF_p(T)$. 
\par
Now take an equation from the table, consider it as a curve over 
$\FF_q(T)$ and assume there is another curve $E'$ over $\FF_q(T)$ 
that becomes isomorphic to $E$ over $\overline{\FF_q}(T)$. Then 
$E'$ can only be the unramified quadratic twist of $E$. Hence 
$E'$ has non-split multiplicative reduction at $\infty$. This 
shows that the table is complete.
\par
A curve from the class $E_1$, having potentially multiplicative 
reduction at $T$, cannot be isogenous to a curve from one of the 
other classes, as those have potentially good reduction at $T$.
For other possible isogenies it suffices if we consider 
$\ell$-isogenies where $\ell$ is a prime. From the theory of Tate 
curves we know that at every pole of $j(E)$ under an $\ell$-isogeny
the pole order will be either multiplied or divided by $\ell$
(regardless of whether the reduction at this place is split or
non-split multiplicative or even additive). This shows that there
cannot be isogenies between different $E_2$, $E_3$ and $E_4$.
\par
Each of the three $\FF_q(T)$-rational $2$-torsion points of 
$(222^*)$ gives rise to a $2$-isogeny to another Frobenius minimal 
elliptic curve in class $E_1$, that is to one of the other three
curves in the box. Similarly, there is a $2$-isogeny between $(12)$
and $(21)$ coming from the $2$-torsion point $(0,0)$.
Finally, the $3$-torsion point $(-9T^2,12T^2(T-1)\sqrt{-3})$ on 
$(13)$ generates an $\FF_q(T)$-rational $3$-isogeny to $(31)$.
\hfill $\Box$
\\ \\
{\bf Remark 2.3.} 
The equations of the isogeny classes $E_1$ and $E_2$ in
Theorem 2.2 also make sense in characteristic $3$, and 
indeed, if $q$ is a power of $3$, there exist only these two 
isogeny classes. Compare [Sch6, Proposition 4.2], where one 
should however replace the second equation by its $(-1)$-twist 
to ensure split multiplicative reduction at $\infty$.
\par
In characteristic $2$ there is only one isogeny class, with 
two Frobenius minimal curves (see [Sch7, Theorem 5.4]). It
corresponds to $E_3$ but, of course, one has to take
equations that are not in short Weierstra\ss\ form.
\\ \\
{\bf Theorem 2.4.} \it
\begin{itemize}
\item[a)] If $q$ is a power of $2$, then the equation 
of the (unique) strong Weil curve over $\FF_q(T)$ with 
conductor $\infty\cdot T^2(T-1)$ is
$$Y^2+XY=X^3+\frac{1}{T^2}X^2+\frac{(T-1)^2}{T^8}.$$
\item[b)] If $q$ is a power of $3$, then there are two 
strong Weil curves over $\FF_q(T)$ with conductor 
$\infty\cdot T^2(T-1)$, namely
$$Y^2 =X^3 +T(T+1)X^2 +T^2 X,$$
which is the curve $(411^*)$ from the table in Theorem 2.2
and represents the isogeny class with supersingular reduction
at the place $T+1$; and
$$Y^2 =X^3 +T^2 X^2 +TX,$$
which is the Frobenius of the curve $(21)$ and represents
the isogeny class with ordinary reduction at $T+1$.
\item[c)] If $char(\FF_q)\ge 7$, then the last equation
in each horizontal box in Theorem 2.2 gives the strong
Weil curve of the corresponding isogeny class.
\par
This also holds in characteristic $5$, except for the
class $E_4$ whose strong Weil curve then is
$$Y^2=X^3+3T^4X-T(T-2)^5,$$
which is the image of the last curve under the Frobenius
isogeny.
\end{itemize}
{\bf Proof. }\rm\\
a) By [Sch7, Theorem 5.4] there are two Frobenius minimal 
curves in this isogeny class, namely 
\begin{eqnarray*}
               & (13):\  & Y^2 +XY=X^3 +\frac{1}{T}X^2 +\frac{(T-1)^3}{T^4}\\
\hbox{\rm and} & (31):\  & Y^2 +XY=X^3 +\frac{1}{T}X^2 +\frac{T-1}{T^4}.
\end{eqnarray*}
The reductions of the second equation modulo $T-a$ with 
$1\neq a\in\FF_q^*$ inherit the rational $3$-torsion point 
$(\frac{1}{T},\frac{1}{T^2})$. By Lemma 1.3 they also have a rational 
$2$-torsion point. But for $q>2$ there exists $a\in\FF_q^*$, 
$a\neq 1$ such that $\frac{1}{a}$ (the coefficient of 
$X^2$) is not of the form $b^2 +b$; so by Lemma 1.3 the 
reduction modulo $T-a$ has no $4$-torsion point.
In the terminology of Corollary 1.2 we thus have $N=6$ (except for 
$q=2$; then there are no linear places of good reduction).
\par
The multiplicative reduction at $T-1$ is non-split if and only if 
$q$ is an odd power of $2$. In this case $6$ divides $2q+2$. All in 
all we obtain $-v_{\infty}(j(E))=6$ (also for $q=2$). So the strong 
Weil curve is the Frobenius of the second equation.
\\ \\
b) We start with the curve $(222^*)$ from the table in Theorem 2.2. 
From the $j$-invariant we see that the reduction at $T$ is supersingular. 
Consequently $3$ does not divide the pole order of the $j$-invariant 
of the strong Weil curve, which therefore has to be Frobenius minimal 
(cf. Corollary 1.2). On the other hand, the $2$-torsion points survive 
every reduction and $2q+2$ is divisible by $4$. So by Theorem 2.2 and 
Corollary 1.2 the strong Weil curve is $(411^*)$.
\par
For the other isogeny class we have to show $-v_{\infty}(j(E))=6$.
We start with the equation $(21)$. The $2$-torsion point $(0,0)$
survives every reduction. We transform to 
$$Y^2 =X^3 +X^2 -\frac{T-1}{T^3}.$$ 
Then, besides the split multiplicative reduction at $T-1$, we 
see by Lemma 1.3 that the reductions at places $T-a$ have 
a $3$-torsion point, and we have to show that at least one has 
no $9$-torsion point over $\FF_q$. For $q=3$ this is clear from 
the Weil bounds. Writing $W$ for $\frac{1}{T}$, the reductions 
we get for $q>3$ are
$$Y^2 =X^3 +X^2 +w^3 -w^2$$
with $1\neq w\in\FF_q^*$. If they all had an $\FF_q$-rational 
$9$-torsion point, then by Lemma 1.3 the elliptic curve 
$U^3 -U=W^3 -W^2$ would have at least $3(q-2)$ rational points 
over $\FF_q$, which is impossible.
\\ \\
c) First we use Corollary 1.2 to show that up to Frobenius the 
strong Weil curve is always the last equation in each horizontal 
box. Obviously all reductions of $(222^*)$ have a full set of
rational $2$-torsion points and $4$ divides $2q+2$. Similarly 
with the $2$-torsion point of $(12)$. The curve $(13)$ has 
a $3$-torsion point 
$$(-9T^2,12T^2(T-1)\sqrt{-3}).$$
If $q\equiv 1\ mod\ 3$, this gives a rational $3$-torsion point 
on the reductions. Moreover, the harmonicity condition at the 
vertex $e_0^1$ of the quotient graph (or equivalently, the fact 
that the $L$-polynomial of our curves are constant $1$) implies
that the multiplicative reduction at $T-1$ is split. 
If $q\equiv 2\ mod\ 3$, then $2q+2$ is divisible by $3$, but the
above point does not give a rational point on the reduced curve
$\widetilde{E}$. However, it shows that the twist of $\widetilde{E}$
has an $\FF_q$-rational $3$-torsion point. Since the number of 
rational points on $\widetilde{E}$ and its twist add up to $2q+2$,
this shows that $\widetilde{E}$ must have (another) rational
$3$-torsion point.
\par
Now we want to show that the strong Weil curve is Frobenius minimal
(at least for $p\ge 7$). Let $A_p(T)$ be the Hasse invariant of $E$.
For our curves $A_p(T)$ has degree at most $p-1$ and is divisible by
$T^2$. By [Si, pp.141/142] the number of $\FF_q$-rational points on 
the reduction of $E$ mod $T-a$ is congruent to $1-A_q(a)$ modulo $p$ 
where 
$$A_q(a)=(A_p(a))^{\frac{q-1}{p-1}}.$$
Now assume that $p$ divides $v_{\infty}(j(E))$. Then $A_q(a)=1$
for all $a\in\FF_q \setminus\{0,1\}$. With the properties mentioned 
above one easily shows that then necessarily $A_q(T)=T^{q-1}$ and
hence $A_p(T)=cT^{p-1}$ with $c\in\FF_p$. This means that for the 
curve over $\FF_p(T)$, i.e. for $q=p$, the values of the corresponding
cycle $\varphi$ on the edges from $e_1^1$ to $e_0^1$ are all congruent
to $c$ modulo $p$. By the harmonicity condition this is only possible
if either $c=0$ (in which case we are done) or if $c=1$. The latter 
case means that all reductions $\widetilde{E}$  at places $T-a$ with 
$a\in\FF_p\setminus\{0,1\}$ have an $\FF_p$-rational $p$-torsion point. 
For the classes $E_1$, $E_2$ and $E_3$ the Weil bound then gives 
$$2p\le \#(\widetilde{E}(\FF_p ))\le p+1+2\sqrt{p},$$ 
which is only possible for $p\le 5$.
Despite first appearance to the contrary, the same argument 
can be made to work for the curve $(11)$. We transform it into
$$Y^2 =X^3 -27X+54(1-\frac{2}{T}).$$
If $p\ge 7$ and $X$ runs through $\FF_p$, then $X^3 -27X$ takes at 
least $3$ different values. So we can find $a\in\FF_p \setminus\{0,1\}$
such that for $T=a$ the cubic polynomial has a zero in $\FF_p$. Thus
there exists a reduction with an $\FF_p$-rational $2$-torsion point.
\par
Finally we deal with characteristic $5$. Then the curve $(222^*)$ has 
Hasse invariant $T^2(T^2 -T+1)$, and the classes $E_2$ and $E_3$ have 
supersingular reduction at $T-2$. Transforming the equation $(11)$ into
$$Y^2 =X^3 +3X+\frac{2}{T}-1$$
we see by Lemma 1.3 that all reductions have a rational $5$-torsion
point and that all elliptic curves with a $5$-torsion point, including
those without a $25$-torsion point, occur among these reductions.
Moreover, the multiplicative reduction at $T-1$ is split. So in this 
case $-v_{\infty}(j(E))=5$.
\hfill $\Box$
\\ \\
We avoid getting lost in trying to write down equations for 
the curves with square-free $\nn\in\FF_q[T]$ of degree $3$. 
\par
In the special case of characteristic $2$ all these curves 
have been explicitly determined in Sections 3 and 4 of [Sch7]. 
In a certain sense the proof in that paper is unnecessarily 
complicated. As an alternative strategy one can see from 
Szpiro's conjecture [PeSz] that the Frobenius minimal ones 
give rational elliptic surfaces, and then one can, as in the 
proof of Theorem 2.2, determine the forms over $\FF_q(T)$ of 
the equations in [La1] and [La2]. This approach was already 
outlined in [Ng, Proposition 4]. As one excuse we mention 
that the complete proof [PeSz] of Szpiro's conjecture in 
characteristic $2$ had not been published yet when [Sch7] 
was written. 
\par
In any case, for any finite field $\FF_q$ and any given
$\nn\in\FF_q[T]$ of degree $3$, writing down all elliptic
curves over $\FF_q(T)$ with conductor $\infty\cdot\nn$ is
(at least in principal) mainly a matter of patience. 
We content ourselves with the following example.
\\ \\
{\bf Example 2.5.} 
Let $\nn\in\FF_q[T]$ be a monic irreducible polynomial of 
degree $3$. Then the necessary and sufficient condition for 
the existence of an elliptic curve $E/\FF_q(T)$ with conductor 
$\infty\cdot\nn$ is that by a translation $T\mapsto T+b$ one 
can bring $\nn$ into the form $T^3-c$. In particular, such 
curves exist if and only if $q\equiv 1\ mod\ 3$.
\par
Indeed, suppose that $E/\FF_q(T)$ is such a curve. We may 
assume that it is Frobenius minimal. Over $\overline{\FF_q}$ 
its equation gives a rational extremal elliptic surface, whose 
$j$-invariant has at least $3$ different poles of the same 
order. Up to affine transformation of $T$ this can only be 
a surface $(3333)$ or $(9111)$ from [La1]. But these surfaces 
don't exist in characteristic $3$. Hence we can find a translation 
$T\mapsto T+b$ that transforms $\nn$ into a polynomial
$T^3+c_1T+c_0$. Since the elliptic curves $(3333)$ and
$(9111)$ have conductor $\infty\cdot(T^3-1)$ we see
that $c_1=0$. But irreducible polynomials of the form
$T^3-c$ exist only if $q\equiv 1\ mod\ 3$.
\par
Conversely, if $q\equiv 1\ mod\ 3$, then $\FF_q^*$ contains 
an element $c$ that is not a third power. Then for example
$$Y^2=X^3-3T(T^3+8c)X-2(T^6-20cT^3-8c^2)$$
has conductor $\infty\cdot(T^3-c)$. We have taken the
$(-1)$-twist of the curve in [Ge1, Table 9.3] to make sure
that the multiplicative reduction at $\infty$ is split.
\\ \\
Even without the explicit knowledge of the curves we 
can make the following statement.
\\ \\
{\bf Theorem 2.6.} \it
Semistable strong Weil curves $E/\FF_q(T)$ of conductor 
$\infty\cdot\nn$ with $deg(\nn)=3$ are Frobenius minimal.
\\ \\
{\bf Proof. }\rm
We even prove that $v_\infty(j(E))$ is not divisible by $p$, 
the characteristic of $\FF_q$. For every $a\in\FF_q$ we have
$\langle\varphi,\delta_a -\delta_\infty\rangle\equiv
\varphi(a)-1\ mod\ q$. Thus, if $p$ divides $v_\infty(j(E))$, 
then all entries of $\varphi$ must be congruent to $1$ mod $p$. 
But this contradicts the condition that the sum over these 
entries has to be $0$.
\hfill $\Box$
\\ \\
In characteristic $3$ we can prove more. By 
[Sch6, Proposition 4.1] the elliptic curves with
square-free $\nn\in\FF_{3^r}[T]$ of degree 
$3$ have two or no linear places of supersingular 
reduction. Calculating modulo $3$ we see 
$\langle\varphi,\varphi\rangle\not\equiv 0\ mod\ 3$.
This shows not only that the strong Weil curve is 
Frobenius minimal, but also that the degree of the 
strong Weil uniformization is not divisible by $3$.
\\ \\
Using additional machinery, Theorem 2.6 can be generalized 
as follows.
\\ \\
{\bf Lemma 2.7.} \it
Semistable strong Weil curves $E/\FF_q(T)$ of conductor 
$\infty\cdot\nn$ where $\nn$ has an irreducible factor 
$\pp$ with $deg(\frac{\nn}{\pp})\le 2$ are Frobenius minimal.
\\ \\
{\bf Proof. }\rm
If $\nn$ is irreducible, this is [Pa1, Theorem 1.2]. More 
generally, in the terminology of [Pa2], if $E$ is not Frobenius 
minimal, then by [Pa2, Theorem 1.1] we have $p\in C(\pp)$,
but $C(\pp)=\emptyset$ since $S(\nn)^{\pp-old}=0$. 
\hfill $\Box$
\\ \\
Comparing Theorems 2.1, 2.4 and 2.6 with the elliptic 
surfaces in [La1] and [La2], or alternatively, with the 
results in [Ito] and [Sch6], we obtain the following fact.
\\ \\
{\bf Theorem 2.8.} \it
If $E/\FF_q(T)$ is a strong Weil curve with conductor
$\infty\cdot\nn$ where $deg(\nn)=3$, then the 
corresponding elliptic surface is a rational surface 
over $\overline{\FF_q}$.
\\ \\
\rm
One might conjecture that the same statement holds for 
$deg(\nn)=4$. This is true at least for $q=2$ (see the table in 
[Sch2]) or if $\nn$ is square-free but does not split completely 
into linear factors. To see the last claim note that for 
$deg(\nn)=4$ in general Szpiro's conjecture [PeSz] still implies 
that Frobenius minimal elliptic curves $E/\FF_q(T)$ with conductor 
$\infty\cdot\nn$ give rational elliptic surfaces over 
$\overline{\FF_q}$, and by Lemma 2.7 those where $\nn$ is 
square-free of degree $4$ but has a non-linear irreducible 
factor are Frobenius minimal.
\par
In any case, the elliptic surfaces corresponding to the strong 
Weil curves with $deg(\nn)=4$ are at least unirational. This also 
implies, by the way, that for these elliptic curves the conjecture 
of Birch and Swinnerton-Dyer holds (see [Shi, Section 3]).
\\

\subsection*{3. Complements}

At least in one case we want to make Theorem 2.6 more 
explicit.
\\ \\
{\bf Example 3.1.}
In characteristic $2$, semistable elliptic curves 
$E/\FF_q(T)$ with conductor $\infty\cdot\nn$ where
$\nn$ is the product of $3$ different linear factors
exist if and only if $q$ is a power of $4$. Moreover, 
for these curves $\nn$  can be transformed to 
$(T-1)(T-s)(T-s^2)$ where $s$ is a primitive third 
root of unity (see Section 3 in [Sch7]). For example, 
$$(3333):\ \ \ Y^2+TXY+Y=X^3+T^3+1$$
is such a curve. Since all its $3$-torsion points are 
rational, one can show by the same arguments as in the 
previous theorems that the corresponding strong Weil 
curve is 
$$(9111):\ \ \ Y^2+TXY+Y=X^3.$$
The numbers $(klmn)$ give the pole orders of the $j$-invariant
at the places $\infty$, $1$, $s$ and $s^2$.
There are $3$ more isogeny classes with this
conductor and their strong Weil curves are
\def\abstand{$\vrule width 0pt height 12.5pt $}
$$ \begin{array}{cl} 
\abstand (5511):\ & Y^2+TXY+Y=X^3+X^2+T\\
\abstand (5115):\ & Y^2+sTXY+Y=X^3+X^2+sT\\
\abstand (5151):\ & Y^2+s^2TXY+Y=X^3+X^2+s^2T\\
\end{array} $$ \\
in the terminology of [Sch7, Theorem 3.2].
\par
If $q=2^{2n+1}$, the curve $(9111)$ has conductor
$\infty\cdot (T-1)(T^2+T+1)$ but in general it will not be 
the strong Weil curve in this situation. For example if 
$q=2$, the strong Weil curve of this class is $(3333)$ 
(see [GeRe, Example 9.7.4] or [Ge2, Example 4.4] or 
[Ge3, Example 9.4]).
\\ \\
{\bf Proposition 3.2.} \it
If $q=4$ and $\nn$ is the product of three different linear
factors, then the curve $X_0(\nn)$ has no $\FF_4(T)$-rational
points except the $8$ cusps.
\\ \\
{\bf Proof. }\rm
After translation we can suppose $\nn=T(T-1)(T-v)$. From
Table 10.2 in [Ge1] we easily calculate that $X_0(T(T-1)(T-v))$
maps with degree $4$ to the corresponding strong Weil curves.
By the previous example, one of these is a transformation of
$(9111)$, which has only $3$ rational points. Thus $X_0(\nn)$ 
has at most $12$ rational points over $\FF_4(T)$. 
\par
But at the same time the number of these points is divisible 
by $8$. This follows from the action of the Atkin-Lehner 
involutions on the rational points. Namely, by the proof 
of Lemma 12 in [Sch4] the fixed 
points of these Atkin-Lehner involutions correspond to 
Drinfeld modules with complex multiplication by orders in 
$\FF_q[\sqrt{T}]$, and by [Sch3, Lemma 4] these Drinfeld 
modules have $j$-invariants that are inseparable over 
$\FF_q(T)$. Thus the fixed points are not rational.
So the Atkin-Lehner involutions (which form a group of 
order $8$) act freely on the $\FF_q(T)$-rational points
of $X_0(\nn)$.
\hfill $\Box$
\\ \\
By analogous arguments one can show that the curve
$X_0(\nn)$ has no $\FF_q(T)$-rational points outside
the cusps for the following values
\begin{eqnarray*}
q=2, & \nn= & T^3,\ \ T^2(T-1),\ \ T(T^2+T+1),\\
q=3, & \nn= & T(T-1)(T+1),\ \ T^2(T-1).
\end{eqnarray*}
But for bigger $q$ the attempt is too weak, simply because
the degree of the strong Weil uniformization grows with
$q$. More precisely:
\\ \\
{\bf Lemma 3.3.} \it
Let $\pi: X_0(\nn)\to E$ be the uniformization of a strong
Weil curve $E$ over $\FF_q(T)$ of conductor $\infty\cdot\nn$ 
with $deg(\nn)=3$. Then
$$\frac{q}{2}\le deg(\pi)\le 
-v_\infty(j(E))deg(\pi)\le 4q^2+q+1.$$
{\bf Proof. }\rm
Let $W_\nn$ be the full Atkin-Lehner involution of 
$X_0(\nn)$. It is well known that in the case $deg(\nn)=3$
the quotient curve $W_\nn\!\setminus\! X_0(\nn)$ is rational
(cf. [Ge1] or [Sch4]).
Applying Castelnuovo's inequality (see for example 
[Sti, Theorem III.10.3]) to $\pi$ and the canonical map 
$\kappa: X_0(\nn)\to W_\nn\!\setminus\! X_0(\nn)$, we obtain
\begin{eqnarray*}
g(X_0(\nn)) & \le & deg(\pi)g(E)
+deg(\kappa)g(W_\nn\!\setminus\! X_0(\nn))
+(deg(\pi)-1)(deg(\kappa)-1)\\
 & = & 2 deg(\pi)-1.
\end{eqnarray*}
For $deg(\nn)=3$ it is also well known that $g(X_0(\nn))=q$
if $\nn$ is square-free, and $g(X_0(\nn))=q-1$ otherwise.
\par
The upper bound comes from estimating 
$\langle\varphi,\varphi\rangle$.
\hfill $\Box$
\\

\subsection*{\hspace*{9em} Acknowledgements}
This paper was written while I was holding a visiting position at the 
National Center for Theoretical Sciences (NCTS) in Hsinchu, Taiwan. 
I want to thank all the people there for their support.

\subsection*{\hspace*{10.5em} References}
\begin{itemize}
\item[{[Ge1]}] E.-U. Gekeler: Automorphe Formen \"uber 
$\FF_q(T)$ mit kleinem F\"uhrer, \it 
Abh. Math. Sem. Univ. Hamburg \bf 55 \rm (1985), 111-146 
\item[{[Ge2]}] E.-U. Gekeler: Analytical Construction of 
Weil Curves over Function Fields, \it 
J. Th\'eor. Nombres Bordeaux \bf 7 \rm (1995), 27-49
\item[{[Ge3]}] E.-U. Gekeler: Jacquet-Langlands 
theory over $K$ and relations with elliptic curves, 
\it in: Drinfeld Modules, Modular Schemes and Applications, 
\rm Proceedings of a workshop at Alden Biesen, 
September 9-14, 1996, (E.-U. Gekeler, M. van der Put, 
M. Reversat, J. Van Geel, eds.), World Scientific, Singapore, 
1997, pp. 224-257
\item[{[Ge4]}] E.-U. Gekeler: Highly ramified pencils 
of elliptic curves in characteristic two, 
\it Duke Math. J. \bf 89 \rm (1997), 95-107
\item[{[Ge5]}] E.-U. Gekeler: Local and global 
ramification properties of elliptic curves in characteristics 
two and three, \it in: Algorithmic Algebra and Number Theory, 
\rm (B. H. Matzat, G.-M. Greuel, G. Hi\ss ,\ eds.), Springer,
Berlin-Heidelberg-New York, 1998, pp. 49-64
\item[{[GeNo]}] E.-U. Gekeler and  U. Nonnengardt: 
Fundamental domains of some arithmetic groups over function
fields, \it Internat. J. Math. \bf 6 \rm (1995), 689-708
\item[{[GeRe]}] E.-U. Gekeler and M. Reversat: 
Jacobians of Drinfeld Modular Curves, 
\it J. Reine Angew. Math. \bf 476 \rm (1996), 27-93
\item[{[Ito]}] H. Ito: On unirationality of extremal elliptic 
surfaces, \it Math. Annalen \bf 310 \rm (1998), 717-733
\item[{[La1]}] W. Lang: Extremal rational elliptic 
surfaces in characteristic $p$. I: Beauville surfaces, 
\it Math. Z. \bf 207 \rm (1991), 429-438
\item[{[La2]}] W. Lang: Extremal rational elliptic surfaces 
in characteristic $p$. II: Surfaces with three or fewer 
singular fibres, \it Ark. Mat. \bf 32 \rm (1994), 423-448
\item[{[Lei]}] R. Leitl: \it Elliptische Kurven \"uber 
$\FF_q(T)$ mit kleinem F\"uhrer, \rm Diplomarbeit, 
Saarbr\"ucken 1995
\item[{[Ng]}] K. V. Nguyen: On families of curves over 
$\PP^1$ with small number of singular fibres, \it
C. R. Acad. Sci. Paris \bf 326 \rm (1998), 459-463
\item[{[Pa1]}] M. Papikian: Pesenti-Szpiro inequality for
optimal elliptic curves, \it J. Number Theory \bf 114 \rm
(2005), 361-393
\item[{[Pa2]}] M. Papikian: Abelian subvarieties of 
Drinfeld Jacobians and congruences modulo the characteristic,
\it Math. Annalen \bf 337 \rm (2007), 139-157
\item[{[PeSz]}] J. Pesenti and L. Szpiro: In\'egalit\'e du 
discriminant pour les pinceaux elliptiques \`a r\'eductions 
quelconques, \it Compositio Math. \bf 120 \rm (2000), 83-117
\item[{[R\"u]}] H. G. R\"uck: A note on elliptic curves over 
finite fields, 
\it Math. Comp. \bf 49 no. 179 \rm (1987), 301-304
\item[{[Sch1]}] A. Schweizer: \it Zur Arithmetik der
Drinfeld'schen Modulkurven $X_0(\nn)$, \rm Dissertation, 
Saarbr\"ucken 1996
\item[{[Sch2]}] A. Schweizer: Modular automorphisms of the 
Drinfeld modular curves $X_0(\nn)$, \it Collect. Math. \bf
48 \rm (1997), 209-216
\item[{[Sch3]}] A. Schweizer: On singular and supersingular 
invariants of Drinfeld modules, 
\it Ann. Fac. Sci. Toulouse Math. \bf 6 \rm (1997), 319-334
\item[{[Sch4]}] A. Schweizer: Hyperelliptic Drinfeld Modular 
Curves, \it in: Drinfeld Modules, Modular Schemes and 
Applications, \rm Proceedings of a workshop at Alden Biesen, 
September 9-14, 1996, (E.-U. Gekeler, M. van der Put, 
M. Reversat, J. Van Geel, eds.), World Scientific, Singapore, 
1997, pp. 330-343
\item[{[Sch5]}] A. Schweizer: On elliptic curves in 
characteristic $2$ with wild additive reduction, \it Acta Arith.
\bf 91 \rm (1999), 171-180 
\item[{[Sch6]}] A. Schweizer: Extremal elliptic surfaces 
in characteristic $2$ and $3$, \it Manu\-scripta Math. \bf 102 
\rm (2000), 505-521
\item[{[Sch7]}] A. Schweizer: On elliptic curves over
function fields of characteristic two, \it J. Number Theory
\bf 87 \rm (2001), 31-53
\item[{[Shi]}] T. Shioda: Some remarks on elliptic curves over 
function fields, \it Ast\'erisque \bf 209 \rm (1992), 99-114
\item[{[Si]}] J. H. Silverman: \it The Arithmetic of Elliptic
Curves, \rm Springer GTM 106, Berlin-Heidelberg-New York, 1986
\item[{[Sti]}] H. Stichtenoth: \it Algebraic Function Fields 
and Codes, \rm Springer Universitext, Berlin-Heidelberg-New York,
1993
\item[{[Vo]}] J.F. Voloch, Explicit $p$-descent for  elliptic 
curves in characteristic $p$, \it Compositio Math. \bf 74 \rm 
(1990), 247-258
\end{itemize} 

\end{document}